\documentclass[a4paper,fleqn,10pt]{amsart}

\usepackage[english]{babel}

\usepackage{amssymb}
\usepackage[sort]{natbib}
\usepackage{mathrsfs}
\usepackage{mathscinet}
\usepackage{fullpage}
\usepackage[colorlinks,linkcolor=blue,citecolor=blue]{hyperref}
\usepackage{kamil}

\begin{document}
\bibliographystyle{plainnat}
\setcitestyle{numbers}

\title{On the infimum attained by the reflected fractional Brownian motion}

\author{K.\ D\polhk{e}bicki}
\address{Instytut Matematyczny, University of
Wroc\l aw, pl.\ Grunwaldzki 2/4, 50-384 Wroc\l aw, Poland.}
\email{Krzysztof.Debicki@math.uni.wroc.pl}

\author{K.M.\ Kosi\'nski}
\address{Korteweg-de Vries Institute for Mathematics,
University of Amsterdam, the Netherlands; E{\sc urandom},
Eindhoven University of Technology}
\email{Kosinski@eurandom.tue.nl}
\thanks{KK was supported by NWO grant 613.000.701.}



\date{\today}

\begin{abstract}

Let
$\{B_H(t):t\ge 0\}$ be a fractional Brownian motion with Hurst parameter
$H\in(\half,1)$.
For the storage process
$Q_{B_H}(t)=\sup_{-\infty\le s\le t} \left(B_H(t)-B_H(s)-c(t-s)\right)$ we show that,
for any $T(u)>0$ such that
$T(u)=o(u^\frac{2H-1}{H})$,
$$\mathbb P (\inf_{s\in[0,T(u)]} Q_{B_H}(s)>u)\sim\mathbb P(Q_{B_H}(0)>u),$$
as $u\to\infty$.
This finding, known in the literature as the {\it strong Piterbarg property},
goes in line with previously observed
properties of storage processes with self-similar and infinitely divisible input without
Gaussian component.
\end{abstract}

\maketitle

\section{Introduction}
\label{sec:intro}

The analysis of distributional properties of {\it reflected} stochastic processes
is continuously motivated both by theory- and applied-oriented open problems in probability
theory.
In this paper we analyze the asymptotic properties of tail distribution of infimum of
an important class of such processes, that naturally appear in models of storage (queueing)
systems and, by duality to ruin problems,
gained broad interest also in problems arising in
finance and insurance risk; see, e.g., \citep{Norros94,Piterbarg01,Asmussen03,Asmussen10}
or a novel work \citep{HJP13}.

Consider a {\it fluid queue} with infinite buffer capacity,
service rate $c>0$ and the total inflow by time $t$
modeled by a stochastic process with stationary increments $X=\{X(t):t\in\rr\}$.
Following \citet{Reich58}, the stationary {\it storage process} that
describes the stationary buffer content process, has the following representation
\[
Q_{X}(t)=\sup_{-\infty\le s\le t} \left(X(t)-X(s)-c(t-s)\right).
\]
There is a strong motivation for modeling the input process $X$ by a fractional Brownian motion (fBm) $B_H=\{B_H(t):t\in\rr\}$ with
$H>1/2$, i.e.,
a centered Gaussian process with stationary increments, continuous sample paths a.s., and variance function $\sigma^2_{B_H}(t)=t^{2H}$.
On one hand, such structural properties of fBm as \textit{self-similarity} and
\textit{long range dependence}, have been statistically confirmed in
data analysis of many real traffic processes in modern data-transfer networks. On the other hand, in \citep{Taqqu97,STE} it was proven that appropriately
scaled aggregation of large number of (integrated) On-Off input processes with regularly varying tail distribution of successive On-times, converges to an
fBm with $H>1/2$.

The importance of fBm storage processes resulted in
a vast interest of analysis of
the process $Q_{B_H}$. In particular finding the properties
of finite-dimensional
(or at least 1-dimensional) distributions of $Q_{B_H}$ has been a long standing goal;
see \citep{Norros94,Piterbarg01}. The stationarity of increments of $B_H$ implies the stationarity of the process
$Q_{B_H}$, so that, for any fixed $t$, the random variable $Q_{B_H}(t)$ has the same
distribution as $Q_{B_H}(0)$. Nevertheless, apart from the Brownian case $H=\half$,
the exact distribution of $Q_{B_H}(0)$ is not known. Therefore, one usually resorts to
the exact asymptotics of $\prob{Q_{B_H}(0)>u}$, as $u\toi$.
These have been found for the full range of parameter $H\in(0,1)$ in \citep{Husler99},
leading to,
\begin{equation}
\label{eq:asymptotics}
\prob{Q_{B_H}(0)>u}
\sim \sqrt\pi a^\frac{1}{2H}b^{-\half}\mathcal H^{\sup}_{B_H} (Au^{1-H})^{\frac{1-H}{H}}\Psi(Au^{1-H}),\as u,
\end{equation}
where the constants $a,b$ and $A$ can be given explicitly
(see Section 4),
$\mathcal H^{\sup}_{B_H}$ is the so-called \textit{Pickands constant}, and $\Psi(u)$ denotes the right tail of the standard normal distribution.

\citet{Piterbarg01} considered the supremum of the process $Q_{B_H}$ on the interval $[0,T]$ and found the exact asymptotics of
\[
\prob{\sup_{t\in[0,T]}Q_{B_H}(t)>u},\as u,
\]
for the whole range of the parameter $H$. By comparing them with \eqref{eq:asymptotics},
he observed a remarkable property that, for $H>\half$, and any positive function $T=T(u)$ such that $T(u)=o(u^\frac{2H-1}{H})$,
\begin{equation}
\label{eq:Piterbarg}
\prob{\sup_{t\in[0,T]}Q_{B_H}(t)>u}\sim\prob{Q_{B_H}(0)>u},\as u.
\end{equation}
This property is nowadays referred to as the \textit{generalized Piterbarg property};
see \citep{Albin04}. As a corollary from \eqref{eq:Piterbarg} one easily gets that
for any fixed $n>0$ and $t_1,\ldots,t_n\in[0,T]$, with $u\toi$,
\[
\prob{\min_{i=1,\ldots,n}Q_{B_H}(t_i)>u\Big|\sup_{t\in[0,T]}Q_{B_H}(t)>u}
\ge
1-\sum_{i=1}^n\left(1-\frac{\prob{Q_{B_H}(t_i)>u}}{\prob{\sup_{t\in[0,T]}Q_{B_H}(t)>u}}\right)\to 1.
\]
This leads to the natural question, whether the minimum over finite number of points can be substituted with the infimum functional, which then leads to
\begin{equation}
\label{eq:strongPiterbarg}
\prob{\inf_{t\in[0,T]}Q_{B_H}(t)>u}\sim
\prob{\sup_{t\in[0,T]}Q_{B_H}(t)>u},\as u.
\end{equation}
This property shall be referred to as the \textit{strong Piterbarg property}.

The above terminology has been coined by \citet{Albin04}, who, motivated by \citep{Piterbarg01}, considered the case when the input process $X$ belongs to
the class of self-similar infinitely divisible stochastic processes with no Gaussian component. They provide general conditions under which
\eqref{eq:Piterbarg} and \eqref{eq:strongPiterbarg} hold with $Q_{X}$ instead of $Q_{B_H}$. The approach in \citep{Albin04} is based on the assumption
that the L\'evy measure associated with $X$ has heavy tails, which combined with the absence of a Gaussian component allows for more direct and less
delicate methods to be employed. It is the light-tailed nature of the Gaussian distribution that renders the problem of the asymptotics of suprema of
Gaussian processes hard. Furthermore, infima of Gaussian processes (apart perhaps from the Brownian case) have not been considered systematically. On the
high level, the problem stems from the fact that an infimum is, by definition, an intersection of events. If the number of events grows to infinity, then
the intersection is much harder to handle than, for instance, the sum of events (which corresponds to the supremum).

In this paper we derive exact asymptotics of
\begin{equation}
\label{eq:asymptinf}
\prob{\inf_{t\in[0,T]} Q_{B_H}(t)>u}, \as u.
\end{equation}
and
prove that the strong Piterbarg property \eqref{eq:strongPiterbarg} holds for the same
range of functions $T(u)$ as in the generalized Piterbrag property \eqref{eq:Piterbarg}, i.e., $T(u)=o(u^\frac{2H-1}{H})$, $H>\half$. The idea of the
proof is based on finding  the exact asymptotics of
\begin{equation}
\label{eq:dum}
\prob{\Phi(X_u)>u},\as u,
\end{equation}
for a broad class of functionals $\Phi:C(T)\to\rr$ acting on the space $C(T)$ of
continuous functions on compacts $T\subset\rr^d_+$, $d\ge 1$, and a broad class of Gaussian fields $X_u=\{X_u(\bt):\bt\in\rr^d_+\}$. The connection
between \eqref{eq:asymptinf} and \eqref{eq:dum} can be seen by setting $d=1$, $\Phi(f)=\inf_{t\in[0,1]}f(t)$ and $X_u(t)=Q_{B_H}(T(u)t)$, although the
relation is far from straight forward since $Q_{B_H}$ is not Gaussian.

{\it Structure of the paper:}
The exact asymptotics of \eqref{eq:dum} are given in \autoref{lem:Pickands}
(see \autoref{sec:Pickands}),
which is the first contribution of this paper.
Interestingly, the asymptotics of \eqref{eq:dum} involve a new type of constants of the form
\[
\mathcal H^{\Phi}_\eta (T)=\ee\exp(\Phi(\sqrt 2 \eta(\cdot)-\sigma_\eta^2(\cdot))),
\]
where $\eta$ is a Gaussian random field with variance function $\sigma^2_\eta$.
These new constants extend the notion of the classical Pickands'
constants
$\mathcal H^{\sup}_{B_H}(S)=\ee\exp(\sup_{t\in[0,S]}(\sqrt 2 B_H(t)-t^{2H}))$, $S>0$,
dating back to \citet{Pickands69a}.
Recall that
$\mathcal H_{B_H}^{\sup}=\lim_{S\toi}\mathcal H_{B_H}^{\sup}([0,S])/S$
in \eqref{eq:asymptotics}.
In \autoref{thm:strong} (\autoref{sec:SPP}) we give the strong Piterbarg property,
which is the second contribution of this paper.
More precisely, we show that \eqref{eq:strongPiterbarg} holds for $H>\half$
and $T(u)=o(u^\frac{2H-1}{H})$, i.e., the same order of functions for which \eqref{eq:Piterbarg} holds.
In \autoref{sec:proof1} and \autoref{sec:proof2}
we give the proofs of our main results.

\section{Notation}
Before we begin, let us set the notation that will be used throughout the paper.
By $B_H=\{B_H(t):t\in\rr\}$ we denote the fBm with Hurst
parameter $H\in(0,1)$, that is, a Gaussian process with zero mean and covariance function given by
\[
\Cov(B_H(t),B_H(s))=\half\left(|t|^{2H}+|s|^{2H}-|t-s|^{2H}\right).
\]
Let $\Psi$ be the right tail of the standard normal distribution.
Recall that
\begin{equation}
\label{eq:Psi}
\Psi(u)=\frac{1}{u\sqrt{2\pi}}\exp\left(-\frac{u^2}{2}\right)\left(1+O(u^{-2})\right),\as u.
\end{equation}
For any vector $\bt\in\rr^d$, $d\ge 1$, we denote $\bt =(t_1,\ldots, t_n)$.
By $\eta=\{\eta(\bt):\bt\in\rr^d_+\}$ we denote a centered Gaussian field, with almost surely continuous sample paths, $\eta(\0)=0$ and variance function
$\sigma_\eta^2(\bt)=\Var(\eta(\bt))$.  Let us introduce the following condition:

\vspace{2mm}

\noindent \textbf{E1:}\ \ \
$\ee \left(\eta(\bt_1)-\eta(\bt_2)\right)^2\le G\|\bt_1-\bt_2\|^\gamma$, for some $\gamma,G>0$ and every $\bt_1,\bt_2\in\rr^d_+$.

\vspace{2mm}

Condition \textbf{E1} is a standard regularity requirement; see, e.g., \citep{Piterbarg96}.
Now let $\Phi:C(T)\to\rr$ be a functional acting on
$C(T)$, the space
of continuous functions on compacts $T\subset\rr^d_+$, $d\ge1 $.
Assume that:
\vspace{2mm}

\noindent \textbf{F1:}\ \ \
$|\Phi(f)|\le \sup_{\bt\in T}f(\bt)$,

\vspace{2mm}

\noindent \textbf{F2:}\ \ \
$\Phi(af+b)=a\Phi(f)+b$, for every $a,b>0$.

\vspace{2mm}

\noindent For $\Phi$ satisfying \textbf{F1} we
define a constant $\mathcal H_{\eta}^{\Phi}(T)$ via
\[
\mathcal H_{\eta}^{\Phi}(T)=
\ee \exp\left(
\Phi\left(\sqrt{2}\eta(\cdot)-\sigma_{\eta}^2(\cdot)\right)
\right).
\]
Note that the dependence on $T$ is implicit via $\Phi:C(T)\to\rr$. To see that the above constant is well defined, notice that due to \textbf{F1},
$\prob{\Phi\left(\sqrt{2}\eta(\cdot)-\sigma_{\eta}^2(\cdot)\right)>u}\le
\prob{\sup_{\bt\in T}\eta(\bt)>u/\sqrt 2}$. Now since $\eta$ is continuous, then it has bounded sample paths a.s. and $\sigma^2_\eta=\sup_{\bt\in
T}\sigma_\eta^2(\bt)<\infty$.  Let $m=\ee\sup_{\bt\in T}\eta(\bt)$. Borell's inequality; see, e.g., \citep{Adler90}, implies that for $x>m$,
$
\prob{\sup_{\bt\in T}\eta(\bt)>x}\le 2\exp\left(-(x-m)^2/(2\sigma_\eta^2)\right)
$
and, as a consequence,
$
\mathcal H_{\eta}^{\Phi}(T)=\int_{-\infty}^\infty e^x \prob{\Phi\left(\sqrt{2}\eta(\cdot)-\sigma_{\eta}^2(\cdot)\right)>x}\D x<\infty.
$

\section{Generalized Pickands' lemma}
\label{sec:Pickands}
In this section we present a lemma that shall play a crucial role in proving the strong Piterbarg property in the remaining part of the paper.

Let us recall that
the original Pickands' lemma \citep{Pickands69,Pickands69a} concerns with
a stationary Gaussian process $X$ with zero mean and covariance function $r(t)$
satisfying $r(t)=1-|t|^{2H}+o(|t|^{2H})$, as $t\to0$, for some $H\in(0,1)$, and $r(t)<1$ for all $t>0$. Its conclusion states that, for any $S>0$,
\begin{equation}
\label{eq:d1}
\prob{\sup_{t\in[0,S]}X_u(t)>u}\sim \mathcal H^{\sup}_{B_{H}}([0,S])\Psi(u),\as u,
\end{equation}
where $X_u(t)=X(tu^{-1/H})$.
Pickands' lemma has been generalized in various ways, capturing both nonstationarity of $X$ and extension to Gaussian fields; see, e.g., \citet{Piterbarg96}.
\citet{Debicki02} presented an extension covering broader local
covariance structures, than
satisfying $\Cov(X(s),X(t))=1-|s-t|^\alpha+o(|s-t|^\alpha)$
as $s-t\to0$, for some $\alpha\in(0,2]$.
Among others, notable extensions have been recently
considered in \cite{DeT11}.

In the following lemma we present a version of  Pickands' lemma
that captures the new constant
$\mathcal H_\eta^{\Phi}(T)$
introduced in the previous section.

\begin{lem}[Generalized Pickands' lemma]
\label{lem:Pickands}
For any $u>0$, let $X_u=\{X_u(\bt):\bt\in\rr^d_+\}$ be a centered Gaussian field with a constant variance equal to one. Let the correlation function
$r_u(\bt_1;\bt_2)=\Corr (X_u(\bt_1), X_u(\bt_2))$ satisfy
\begin{equation}
\label{eq:correlation}
\lim_{u\toi}\sup_{\bt_1,\bt_2\in T}\left|\frac{f^2(u)\left(1-r_u(\bt_1;\bt_2)\right)}{\Var(\eta(\bt_1)-\eta(\bt_2))}-1\right|=0,
\end{equation}
for some compact set $T\subset\rr^d_+$, some function $f(u)\toi$, as $u\toi$,
and $\eta$ satisfying \textbf{E1}. Let $\Phi:C(T)\to\rr$ be a functional satisfying \textbf{F1-F2}. Then, for any function $n(u)$ such that $n(u)\sim
f(u)$,
\[
\prob{\Phi(X_u)>n(u)}
\sim
\mathcal H_{\eta}^{\Phi}(T)\Psi(n(u)),\as u.
\]
\end{lem}
\begin{remark}
Conditions similar to assumption \eqref{eq:correlation}
have been introduced in, among others, \citep{Debicki02,Husler04,Debicki08,DeT11} as a standard way of capturing nonstationarity. The shape of
\autoref{lem:Pickands} is tailored to the needs of the next section, where asymptotics of tail distribution of $\inf\sup$ functionals of Gaussian
processes are
analyzed.
Various further extensions of \autoref{lem:Pickands} can be thought of along the lines
of already existing extensions of the classical Pickands' lemma,
especially in the direction
allowing nonconstant variance function of the family $(X_u)$, as in \citet{Piterbarg78}
or \citet{HJP13}.
\end{remark}

\begin{example}
\label{rem:specialcorr}
Assume that $X=\{X(\bt):\bt\in\rr^d_+\}$ is a centered Gaussian field with unit variance and correlation function satisfying
\[
r(\bt_1;\bt_2)=1-\sum_{i=1}^d a_i|t_{1,i}-t_{2,i}|^{2H_i}+
o\left(\sum_{i=1}^d|t_{1,i}-t_{2,i}|^{2H_i}\right),\text{ as } \sum_{i=1}^d|t_{1,i}-t_{2,i}|\to 0,
\]
for some $H_i\in(0,1)$, $a_i>0$, $i=1,\ldots,d$.
Define a new field $X_u=\{X_u(\bt):\bt\in\rr^d_+\}$ via
$
X_u(\bt)=X(t_1u^{-\frac{1}{H_1}},\ldots,t_du^{-\frac{1}{H_d}}).
$
For any compact set $T\subset\rr^d_+$, the process $X_u$ satisfies \eqref{eq:correlation} with $f(u)=u$ and
$
\bs\eta(\bt)=\sum_{i=1}^d B^i_{H_i}\left(a^\frac{1}{2H_i}t_i\right),
$
where $B^i_{H_i}$ constitute independent fBm's with Hurst parameters $H_i$. Hence the conclusion of \autoref{lem:Pickands} holds for any functional $\Phi$
on $C(T)$ satisfying \textbf{F1-F2}.
In the following section we shall encounter this example in the setting of $d=2$, $H_1=H_2$, $a_1=a_2$ and
$\Phi(f)=\inf_{t_1\in[0,\lambda_1]}\sup_{t_2\in[0,\lambda_2]}f(\bt)$, for some $\lambda_1,\lambda_2>0$. In this case, with $H=H_1$ and $a=a_1$, for any
function $n(u)\sim u$,
\[
\prob{\inf_{t_1\in[0,\lambda_1]}\sup_{t_2\in[0,\lambda_2]}X_u(\bt)>n(u)}
=
\mathcal H_{B_H}^{\inf}([0,a^\frac{1}{2H}\lambda_1])
\mathcal H_{B_H}^{\sup}([0,a^\frac{1}{2H}\lambda_2])\Psi(n(u)),\as u.
\]
\end{example}

\section{Strong Piterbarg property}
\label{sec:SPP}
In this section we present the main result of this paper.
Let us first recall the definition of the storage process $Q_{B_H}$ with service
rate $c>0$ and input $B_H$,
\[
Q_{B_H}(t)=\sup_{-\infty\le s\le t} \left(B_H(t)-B_H(s)-c(t-s)\right).
\]
Let us define the following constants: $a=\frac{1}{2}\tau_0^{-2H}$, $b=\frac{B}{2A}$,
$A=\frac{1}{1-H}\tau_0^{-H}$, $B=H\tau_0^{-H-2}$, $\tau_0=\frac{H}{c(1-H)}$, see
\eqref{eq:asymptotics}. Finally, let
\begin{equation}
\label{eq:Pickands}
\mathcal H^{\sup}_{B_H}=\lim_{S\toi}\frac{\mathcal H_{B_H}^{\sup}([0,S])}{S}
\end{equation}
be the classical Pickands's constant. Now we are in position to state our main result.
\begin{theorem}[Strong Piterbarg property]
\label{thm:strong}
For $H>\half$ and any $T(u)>0$, such that $T(u)=o(u^{\frac{2H-1}{H}})$,
\[
\prob{\inf_{t\in[0,T(u)]}Q_{B_H}(t)>u}
\sim
\sqrt\pi a^\frac{1}{2H}b^{-\half}\mathcal H^{\sup}_{B_H}\cdot (Au^{1-H})^{\frac{1-H}{H}}\Psi(Au^{1-H}),\as u.
\]
In particular,
\[
\prob{\inf_{t\in[0,T(u)]}Q_{B_H}(t)>u}
\sim
\prob{Q_{B_H}(0)>u}
\sim\prob{\sup_{t\in[0,T(u)]}Q_{B_H}(t)>u}, \as u.
\]
\end{theorem}
\begin{remark}
The asymptotics of $\prob{Q_{B_H}(0)>u}$ were found in \citep[Theorem 1]{Husler99}; cf. \eqref{eq:asymptotics}.
The asymptotic equivalence between the tail decay of the supremum functional and the value of $Q_{B_H}$ at $0$ was proven in \citep[Theorem
5]{Piterbarg01} and is called the {\it Piterbarg property}, as mentioned in the introduction; cf. \eqref{eq:Piterbarg}.
Note that the formula in \citet[Theorem 5]{Piterbarg01} should have $a^{\frac{1}{2H}}$ as cited here instead of $a^{\frac{1}{H}}$.
\end{remark}
\begin{remark}
The case of Brownian motion, that is $H=\half$, has been treated in \citep[Theorem 3]{Debicki12a}. The authors found the exact distribution of the infimum
of $Q_{B_\half}$ attained on any interval of the form $[0,S]$, $S>0$,
\[
\prob{\inf_{t\in[0,S]}Q_{B_\half}(t)>u}=\prob{Q_{B_\half}(0)>u}
\left(
2(1+S)\Psi(\sqrt S)-\sqrt{\frac{2S}{\pi}}\exp\left(-\frac{S}{2}\right)
\right).
\]
Let us recall that $Q_{B_\half}(0)$ has $\half$-exponential distribution.
On the other hand, \citep[Theorem 6]{Piterbarg01}, gives
(note that the original formula in \citep{Piterbarg01} has a misprint)
\[
\prob{\sup_{t\in[0,S]}Q_{B_\half}(t)>u}\sim\prob{Q_{B_\half}(0)>u}2\sqrt\pi\mathcal H_{B_\half}^{\sup}([0,2S]), \as u.
\]
Therefore, we see that the strong Piterbarg property does not hold in the case of $H=\half$.
\end{remark}

\begin{remark}
One can envision that the strong Piterbarg property can be applied to functionals $\Phi:C([0,T])\to\rr$ of $Q_{B_{H}}$ that can be majorized, up to the
same magnitude, by the infimum and supremum functionals. A simple example is the integral functional. \autoref{thm:strong} yields, for every $H>\half$,
\[
\prob{\int_{0}^{T(u)} Q_{B_H}(t)\D t>u}\sim \prob{Q_{B_H}(0)>\frac{u}{T(u)}},\as u,
\]
for every function $T(u)>0$ such that $T(u)=o(u^{\frac{2H-1}{3H-1}})$. The problem of the area under the graph of the storage process fed by the Brownian
motion, i.e., the case when $H=\half$, has been considered in \citep{Arendarczyk13}.
\end{remark}

\section{Proof of \autoref{lem:Pickands}}
\label{sec:proof1}
The general idea behind the proof follows the one in \citet[Lemma D.2]{Piterbarg96}.
For any $u>0$,
\begin{align*}
\lefteqn{\prob{\Phi(X_u)>n(u)}
=
\frac{1}{\sqrt{2\pi}}
\int_{\rr}\exp\left(-\frac{v^2}{2}\right)
\prob{\Phi(X_u)>n(u)\Big| X_u(\bs 0)=v}\D v}\\
&\sim
\Psi(n(u))
\int_{\rr}\exp\left(w-\frac{w^2}{2n^2(u)}\right)
\prob{\Phi(X_u)>n(u)\Big| X_u(\bs 0)=n(u)-\frac{w}{n(u)}}\D w,
\end{align*}
where we have used the change of variable $v=n(u)-\frac{w}{n(u)}$.
Let $\zeta_u=\{\zeta_u(\bt):\bt\in T\}$ be a Gaussian field defined via
$
\zeta_u(\bt)=n(u)(X_u(\bt)-n(u))+w.
$
Then, using \textbf{F2}, the last integral can be written as
\[
\int_{\rr}\exp\left(w-\frac{w^2}{2n^2(u)}\right)
\prob{\Phi(\zeta_u)>w
\Big|\zeta_u(\bs 0)=0}\D w
=
\int_{\rr}\exp\left(w-\frac{w^2}{2n^2(u)}\right)
\prob{\Phi(\chi_u)>w}\D w,
\]
where $\chi_u=\{\chi_u(\bt):\bt\in T\}$ is a Gaussian field defined as
$
\chi_u(\bt)\de \zeta_u(\bt)|\zeta_u(\bs 0)=0.
$
For the family of Gaussian distributions that appear inside the integral, for every $\bt\in T$,
\begin{align}
\ee\chi_u(\bt)&=
n(u) \ee\left(X_u(\bt)\left|\right.X_u(\bs 0)=n(u)-\frac{w}{n(u)}\right)-n^2(u)+w\nonumber\\
\label{eq:expectation1}
&=
-n^2(u)(1-r_u(\bt;\bs 0))+ w(1-r_u(\bt;\bs 0)),\\
\ee\chi_u(\bs 0)
&=
\ee\chi_u^2(\bs 0)=0.\nonumber
\end{align}
Furthermore, for any $\bt_1,\bt_2\in T$,
\begin{align*}
\Var\left(\chi_u(\bt_1)-\chi_u(\bt_2)\right)
&=
n^2(u)\left(\Var\left((X_u(\bt_1)-X_u(\bt_2))\Big|X_u(\bs 0)=u-\frac{w}{u}\right)\right)\\
&=
2n^2(u)\left(1-r_u(\bt_1;\bt_2)\right)
-n^2(u)\left(r_u(\bt_1;\bs 0)-r_u(\bt_2;\bs 0)\right)^2.
\end{align*}
Hence from \eqref{eq:correlation} it follows that, as $u\toi$, uniformly on $T$,
\begin{align}
\label{eq:expectation2}
\ee\chi_u(\bt)&\to - \sigma_\eta^2(\bt),\\
\label{eq:variance}
\Var\left(\chi_u(\bt_1)-\chi_u(\bt_2)\right)&\to 2
\Var(\eta(\bt_1)-\eta(\bt_2)).
\end{align}
Thus the finite dimensional distributions of $\chi_u$ converge to the finite
dimensional distributions of
$\tilde\eta=\{\sqrt 2\eta(\bt) -\sigma_\eta^2(\bt):\bt\in T\}$.
Therefore $\chi_u\dto\tilde\eta$ in $C(T)$, as $u\toi$, provided that the family
$\chi=\{\chi_u:u> 0\}$ is tight. For this
let $\chi_u^\circ=\{\chi_u^\circ(\bt):\bt\in T\}$ be a centered Gaussian field defined by $\chi_u^\circ(\bt)=\chi_u(\bt)-\ee\chi_u(\bt)$.
In order to prove tightness of the family $\chi=\{\chi_u:u>0\}$ it suffices to show tightness of the centered family $\chi^\circ=\{\chi^\circ_u:u>0\}$.
Since $\chi_u^\circ(\0)=0$ for all $u>0$, then a straightforward consequence of Straf's criterion for tightness of Gaussian fields, \citep{Straf72},
implies that it suffices to show that for any $\mu,\rho>0$, there exists $\delta\in(0,1)$ and $u_0>0$ such that, for each $\bt_1\in T$ and $u>u_0$,
\begin{equation}
\label{eq:tight}
\prob{\sup_{\|\bt_1-\bt_2\|\le \delta}\left|\chi^\circ_u(\bt_1)-\chi^\circ_u(\bt_2)\right |\ge\mu}
\le
\rho\delta^d,
\end{equation}
where $\|\bt\|=\max\{|t_1|,\ldots,|t_d|\}$. Note that, for sufficiently large $u$,
\[
\ee\left(\chi_u^\circ(\bt_1)-\chi_u^\circ(\bt_2)\right)^2
\le C \Var(\eta(\bt_1)-\eta(\bt_2)),
\]
for all $\bt_1,\bt_2\in T$ and some constant $C>0$. Thus, the assumption \textbf{E1} implies,
\[
\sup_{\|\bt_1-\bt_2\|\le\delta}\Var\left(\chi_u^\circ(\bt_1)-\chi_u^\circ(\bt_2)\right)
\le
CG\delta^{\gamma},
\]
which combined with the application of Borell's inequality gives \eqref{eq:tight}.

Then, the continuous mapping theorem implies
\begin{align}
\label{eq:limit}
\lefteqn{\lim_{u\toi}
\int_{\rr}\exp\left(w-\frac{w^2}{2n^2(u)}\right)
\prob{\Phi(\chi_u)>w}\D w
=
\int_{\rr}
\exp(w)
\prob{\Phi
\left(\sqrt 2\eta(\cdot) -\sigma_\eta^2(\cdot)\right)
>w}\D w}\\
&=
\ee \exp\left(
\Phi\left(\sqrt 2\eta(\cdot) -\sigma_\eta^2(\cdot)\right)
\right)
=\mathcal H_\eta^\Phi(T)\nonumber,
\end{align}
provided we can interchange the limit with the integral in \eqref{eq:limit}.
From \eqref{eq:correlation} it follows that $(1-r_u(\bt;\bs 0))\to 0$ uniformly in $\bt\in T$, therefore \eqref{eq:expectation1}--\eqref{eq:expectation2}
imply that for any $\varepsilon>0$ and sufficiently large $u$,
\[
w_u:=\sup_{\bt\in T}\ee\chi_u^\circ(\bt)\le \varepsilon|w|.
\]
Using \eqref{eq:variance} combined with Sudakov--Fernique's inequality yields, for sufficiently large $u$ and some constant $C>0$,
\[
m_u:=\ee\sup_{\bt\in T}\chi_u^\circ(\bt)
\le C\ee\sup_{\bt\in T}\eta(\bt)=:m.
\]
Furthermore, \eqref{eq:variance} combined with \textbf{E1} implies, for sufficiently large $u$,
\[
\sigma_u^2:=
\sup_{\bt\in T}\Var(\chi^\circ_u(\bt))
\le C\sup_{\bt\in T}\sigma_\eta^2(\bt)\le CG(\diam(T))^{\gamma}.
\]
Now, by \textbf{F1}, Borell's inequality yields, for $|w|(1-\varepsilon)\ge m$,
\begin{align*}
\lefteqn{\prob{\Phi(\chi_u)>w}
\le
\prob{\sup_{\bt\in T}\chi_u^\circ(\bt)>w-w_u}
\le
\prob{\sup_{\bt\in T}\chi_u^\circ(\bt)-m_u>w-\varepsilon |w|-m_u}}\\
&\le
2\exp
\left(
-\frac{(w-\varepsilon |w|-m_u)^2}{2\sigma_u^2}
\right)
\le
2\exp
\left(
-\frac{(w-\varepsilon |w|-m)^2}{2CG(\diam(T))^{\gamma}}
\right).
\end{align*}
Hence the interchange of the limit with the integral in \eqref{eq:limit} follows by the dominated convergence theorem and the limit is finite, that is
$\mathcal H_\eta^\Phi(T)<\infty$.

This completes the proof of Lemma \ref{lem:Pickands}.

\section{Proof of \autoref{thm:strong}}
\label{sec:proof2}
We divide the proof on a number of steps.
Before we proceed, let us make the following observation. The time-reversibility property of fBm implies that (on the process level)
\[
Q_{B_H}(t)\de\sup_{\sigma\ge t}\left(B_H(\sigma)-B_H(t)-c(\sigma-t)\right),
\]
which is the form of $Q_{B_H}$ that we shall use in this section. The relations
of \autoref{sec:reduction} and \autoref{sec:correlation} were derived in \citep{Piterbarg01}.
\subsection{Reduction to a Gaussian field}
\label{sec:reduction}
Using new variables $\tau=(\sigma-t)/u$ and $s=t/u$, for any $T>0$,
\begin{align*}
\lefteqn{\prob{\inf_{t\in[0,T]}Q_{B_H}(t)>u}=
\prob{\inf_{t\in[0,T]}\sup_{\sigma\ge t}\left(B_H(\sigma)-B_H(t)-c(\sigma-t)\right)>u}}\\
&=
\prob{\forall{s\in\left[0,\frac{T}{u}\right]}\exists{\tau\ge 0}:B_H(u(s+\tau))-B_H(s u)>u+cu\tau}\\
&=
\prob{\inf_{s\in[0,Tu^{-1}]}\sup_{\tau\ge 0}\frac{B_H(u(s+\tau))-B_H(s u)}{\tau^Hu^H\nu(\tau)}>u^{1-H}}\\
&=
\prob{\inf_{s\in[0,Tu^{-1}]}\sup_{\tau\ge 0}Z_u(s,\tau)>u^{1-H}},
\end{align*}
where $\nu(\tau)=\tau^{-H}+c\tau^{1-H}$ and $Z_u=\{Z_u(s,\tau):s,\tau\ge 0\}$ is a Gaussian field given by
\[
Z_u(s,\tau)=\frac{B_H(u(s+\tau))-B_H(s u)}{\tau^Hu^H\nu(\tau)}.
\]
The distribution of $Z_u$ does not depend on $u$, hence we deal with $Z=Z_1$. Note that $Z(s,\tau)$ is stationary in $s$, but not in $\tau$.

\subsection{Correlation structure of \texorpdfstring{$Z$}{Z}}
\label{sec:correlation}
The variance $\sigma^2_Z(\tau)$ of $Z(s,\tau)$ equals $\nu^{-2}(\tau)$ and has a single maximum point at $\tau_0=\frac{H}{c(1-H)}$.
Taylor expansion shows that, as $\tau\to\tau_0$,
\begin{equation}
\label{eq:varZ}
\sigma_Z(\tau)=\frac{1}{A}-\frac{B}{2A^2}(\tau-\tau_0)^2+O((\tau-\tau_0)^3),
\end{equation}
where
\begin{align*}
A&=\frac{1}{1-H}\left(\frac{H}{c(1-H)}\right)^{-H}=\nu(\tau_0),\\
B&=H\left(\frac{H}{c(1-H)}\right)^{-H-2}=\nu''(\tau_0).
\end{align*}
Furthermore, denote $a=\frac{1}{2}\tau_0^{-2H}$ and $b=\frac{B}{2A}$. Note that $\tau_0,A,B,a,b$ are the same constants as in \autoref{sec:SPP}.

The correlation function $r(s_1,\tau_1;s_2,\tau_2)$ of $Z$ equals
\begin{align}
\lefteqn{r(s_1,\tau_1;s_2,\tau_2)=\ee Z(s_1,\tau_1) Z(s_2,\tau_2)\nu(\tau_1)\nu(\tau_2)}\nonumber\\
&=
\frac{|s_1-s_2+\tau_1|^{2H}+|s_1-s_2-\tau_2|^{2H}-|s_1-s_2+\tau_1-\tau_2|^{2H}-|s_1-s_2|^{2H}}{2\tau_1^H\tau_2^H}\nonumber\\
&=1-a(1+o(1))\left(|s_1-s_2+\tau_1-\tau_2|^{2H}+|s_1-s_2|^{2H}\right)
\label{eq:corZ}
\end{align}
as $s_1-s_2\to0$, $\tau_1\to\tau_0$, $\tau_2\to\tau_0$.

\subsection{Asymptotic properties of \texorpdfstring{$Z$}{Z}}
In this step we will be concerned with the asymptotic properties of
\begin{equation}
\label{eq:Zasympt}
\prob{\inf_{s\in[0,T]}\sup_{\tau\ge 0} A Z(s,\tau)>u}
\end{equation}
as $u$ grows to infinity. Note that we normalized $Z$ such that now the variance of $AZ(s,\tau)$ equals one at $\tau=\tau_0$ ($Z$ is stationary in $s$).
It follows from \citep[Lemma 1]{Piterbarg01} that there exists a constant $C$ such that, for any $T>0$ and sufficiently large $u$,
\[
\prob{\inf_{s\in[0,T]}\sup_{|\tau - \tau_0|\ge \log u/ u} A Z(s,\tau)>u}
\le
C T u^{2/H}\exp\left(-\half u^2 - b\log^2 u\right).
\]
If we restrict ourselves to the neighborhood
$\{\tau:|\tau-\tau_0|\le\log u/u\}$ of $\tau_0$,
then the following step shows that the probability in \eqref{eq:Zasympt},
with $Z$ restricted to the neighborhood of $\tau_0$,
on the logarithmic scale decays as $-\frac{u^2}{2}$ when $u$ grows large.
Therefore, the neighborhood of $\tau_0$ has the largest contribution
to the asymptotic behavior of \eqref{eq:Zasympt}.
In the following step we present its asymptotic contribution.

\subsection{The asymptotics of the main contributor}
In this step we show that
for any $\lambda>0$, with $\mathcal H_{B_H}^{\sup}$ defined in \eqref{eq:Pickands},
\begin{equation}
\label{eq.inf1}
\liminf_{u\toi}\frac{\prob{\inf_{s\in[0,\lambda u^{-1/H}]}\sup_{|\tau-\tau_0|\le\log u/u} A Z(s,\tau)>u}}
{\sqrt \pi a^\frac{1}{2H}b^{-\half}\mathcal H_{B_H}^{\sup}\mathcal{ H}^{\inf}_{B_H}\left([0,\lambda a^\frac{1}{2H}]\right)u^{\frac{1}{H}-1}\Psi(u)}
\ge 1.
\end{equation}
For the Gaussian field $X(s,\tau)=AZ(s,\tau-s)$,
we have
\[
\prob{\inf_{s\in[0,\lambda u^{-\frac{1}{H}}]}\sup_{|\tau-\tau_0|\le \frac{\log u}{u}}A Z(s,\tau)>u}
\ge
\prob{\inf_{s\in[0,\lambda u^{-\frac{1}{H}}]}\sup_{\tau\in I} X(s,\tau)>u},
\]
for sufficiently large $u$, where
$I:= [\tau_0-\frac{\log u}{2u},\tau_0+\frac{\log u}{2u}]$
(we use that
$I\subset [\tau_0+s- \frac{\log u}{u},\tau_0+s+\frac{\log u}{u}]$ for sufficiently large $u$).
From \eqref{eq:corZ} it follows that the correlation function $r_X$ of $X$ is given by
\[
r_X(s_1,\tau_1;s_2,\tau_2)=1-a(1+o(1))\left(|\tau_1-\tau_2|^{2H}+|s_1-s_2|^{2H}\right)
\]
as $s_1-s_2\to0$, $\tau_1-s_1\to\tau_0$, $\tau_2-s_2\to\tau_0$. Furthermore, \eqref{eq:varZ} implies that the variance function $\sigma_X^2$ of $X$
satisfies
\[
\sigma_X(s,\tau)=1-b(\tau-s-\tau_0)^2+O((\tau-s-\tau_0)^3),
\]
as $\tau-s\to\tau_0$.

Let us divide the interval $[\tau_0- \frac{\log u}{2u},\tau_0+\frac{\log u}{2u}]$ into intervals of length $\gamma u^{-\frac{1}{H}}$ for some fixed
$\gamma>0$,
\begin{align*}
I_k&= [\tau_0 + k\gamma u^{-\frac{1}{H}},\tau_0+(k+1)\gamma u^{-\frac{1}{H}}],\quad k=0,1,2,\ldots,\\
I_{-k}&= [\tau_0-(k+1)\gamma u^{-\frac{1}{H}},\tau_0-k\gamma u^{-\frac{1}{H}}],\quad k=0,1,2,\ldots,
\end{align*}
Notice that,
\begin{align*}
\lefteqn{\prob{\inf_{s\in[0,\lambda u^{-\frac{1}{H}}]}\sup_{\tau\in I} X(s,\tau)>u}}\\
&\ge
\prob{\inf_{s\in[0,\lambda u^{-\frac{1}{H}}]}\max_{k=-[\gamma^{-1} u^{\frac{1}{H}}\frac{\log u}{2u}],\ldots,[\gamma^{-1} u^{\frac{1}{H}}\frac{\log
u}{2u}]}
\sup_{\tau\in I_k} X(s,\tau)>u}\\
&\ge
\prob{\max_{k=-[\gamma^{-1} u^{\frac{1}{H}}\frac{\log u}{2u}],\ldots,[\gamma^{-1} u^{\frac{1}{H}}\frac{\log u}{2u}]}\inf_{s\in[0,\lambda
u^{-\frac{1}{H}}]}\sup_{\tau\in I_k} X(s,\tau)>u}\\
&\ge
2\sum_{k=0}^{[\gamma^{-1}u^{\frac{1}{H}}\frac{\log u}{ 2u}]}
\prob{\inf_{s\in[0,\lambda u^{-\frac{1}{H}}]}\sup_{\tau\in I_k} X(s,\tau)>u}\\
&\quad -
2\sum_{0\le l<k\le [\gamma^{-1}u^{\frac{1}{H}}\frac{\log u}{ 2u}]}
\mathbb P\Bigg(\inf_{s\in[0,\lambda u^{-\frac{1}{H}}]}\sup_{\tau\in I_k} X(s,\tau)>u,
\inf_{s\in[0,\lambda u^{-\frac{1}{H}}]}\sup_{\tau\in I_l} X(s,\tau)>u\Bigg)\\
&\quad -
\mathbb P\Bigg(\inf_{s\in[0,\lambda u^{-\frac{1}{H}}]}\sup_{\tau\in I_{-0}} X(s,\tau)>u,
\inf_{s\in[0,\lambda u^{-\frac{1}{H}}]}\sup_{\tau\in I_0} X(s,\tau)>u\Bigg).
\end{align*}
Now,  for any $\varepsilon>0$, any $s\in[0,\lambda u^{-\frac{1}{H}}]$ and all  $\tau\in I_{\pm k}$, for  sufficiently large $u$,
\[
1-(b+\varepsilon) (k+1)^2\gamma^2 u^{-\frac{2}{H}}\le\sigma_X(s,\tau)\le 1-b(1-\varepsilon) k^2\gamma^2 u^{-\frac{2}{H}}.
\]
Therefore, with $\bar X(s,\tau)= X(s,\tau)/\sigma_X(s,\tau)$,
\[
\prob{\inf_{s\in[0,\lambda u^{-\frac{1}{H}}]}\sup_{\tau\in I_k} X(s,\tau)>u}
\ge
\prob{\inf_{s\in[0,\lambda u^{-\frac{1}{H}}]}\sup_{\tau\in I_k} \bar X(s,\tau)>u_{k+}},
\]
where
\[
u_{k+}=\frac{u}{1-(b+\varepsilon) (k+1)^2\gamma^2 u^{-\frac{2}{H}}}.
\]
Thus by \autoref{rem:specialcorr}, as $u\toi$,
\begin{align*}
\lefteqn{2\sum_{k=0}^{[\gamma^{-1}u^{\frac{1}{H}}\frac{\log u}{ 2u}]}
\prob{\inf_{s\in[0,\lambda u^{-\frac{1}{H}}]}\sup_{\tau\in I_k} X(s,\tau)>u}}\\
&\ge
2(1+o(1))\sum_{k=0}^{[\gamma^{-1}u^{\frac{1}{H}}\frac{\log u}{ u}]}
\mathcal H_{B_H}^{\inf}([0,\lambda a^{\frac{1}{2H}}])
\mathcal H_{B_H}^{\sup}([0,\gamma a^{\frac{1}{2H}}])
\Psi(u_{k_+}).
\end{align*}
Notice that (cf. \eqref{eq:Psi}), as $u\toi$,
\[
\sum_{k=0}^{[\gamma^{-1}u^{\frac{1}{H}}\frac{\log u}{ 2u}]}\Psi(u_{k_+})
\sim
\frac{1}{\sqrt{2\pi}}
\sum_{k=0}^{[\gamma^{-1}u^{\frac{1}{H}}\frac{\log u}{2 u}]}
\frac{1}{u_{k_+}}e^{-\half u^2_{k_+}}.
\]
Furthermore, as $u\toi$,
\begin{align*}
\lefteqn{\frac{1}{\sqrt{2\pi}}\sum_{k=0}^{[\gamma^{-1}u^{\frac{1}{H}}\frac{\log u}{ 2u}]}
\frac{1}{u_{k_+}}e^{-\half u^2_{k_+}}
=
\frac{1}{u\sqrt{2\pi}}
\sum_{k=0}^{[\gamma^{-1}u^{\frac{1}{H}}\frac{\log u}{ 2u}]}
(1-(b+\varepsilon) (k+1)^2\gamma^2 u^{-\frac{2}{H}})}\\
&\quad\times
\exp\left(\frac{-u^2}{2(1-(b+\varepsilon) (k+1)^2\gamma^2 u^{-\frac{2}{H}})^2}\right)\\
&=
\frac{1}{u\sqrt{2\pi}}
\sum_{k=0}^{[\gamma^{-1}u^{\frac{1}{H}}\frac{\log u}{2 u}]}
\exp\left(\frac{-u^2}{2(1-(b+\varepsilon) (k+1)^2\gamma^2 u^{-\frac{2}{H}})^2}\right)(1+o(1))\\
&=
\frac{1}{u\sqrt{2\pi}}
\sum_{k=0}^{[\gamma^{-1}u^{\frac{1}{H}}\frac{\log u}{ 2u}]}
\exp\left(\frac{-u^2(1+(b+\varepsilon) (k+1)^2\gamma^2 u^{-\frac{2}{H}})^2}
                  {2(1-(b+\varepsilon)^2 (k+1)^4\gamma^4 u^{-\frac{4}{H}})^2}\right)(1+o(1))\\
&=
\frac{1}{u\sqrt{2\pi}}
\exp\left(-\frac{u^2}{2}\right)
\sum_{k=0}^{[\gamma^{-1}u^{\frac{1}{H}}\frac{\log u}{ 2u}]}
\exp\left(\frac{-u^2(b+\varepsilon) (k+1)^2\gamma^2 u^{-\frac{2}{H}}}
                  {(1-(b+\varepsilon)^2 (k+1)^4\gamma^4 u^{-\frac{4}{H}})^2}\right)(1+o(1))\\
&=
\Psi(u)
\sum_{k=0}^{[\gamma^{-1}u^{\frac{1}{H}}\frac{\log u}{ 2u}]}
\exp\left(-(b+\varepsilon) (k+1)^2\gamma^2 u^{2-\frac{2}{H}}\right)(1+o(1))\\
&=
\Psi(u) u^{\frac{1}{H}-1}
\sum_{k=0}^{[\gamma^{-1}u^{\frac{1}{H}}\frac{\log u}{ 2u}]}
 u^{1-\frac{1}{H}}\exp\left(-(b+\varepsilon)\gamma^2 \left((k+1) u^{1-\frac{1}{H}}\right)^2\right)(1+o(1))\\
&=
\Psi(u) u^{\frac{1}{H}-1}
\int_0^\infty \exp\left(-(b+\varepsilon)\gamma^2 x^2\right)\D x(1+o(1))
\end{align*}
and
\[
\int_0^\infty \exp\left(-(b+\varepsilon)\gamma^2 x^2\right)\D x =\frac{\sqrt\pi}{2\gamma\sqrt{b+\varepsilon}}.
\]
Combining these estimates we obtain
\begin{align*}
\lefteqn{2\sum_{k=0}^{[\gamma^{-1}u^{\frac{1}{H}}\frac{\log u}{ 2u}]}
\prob{\inf_{s\in[0,\lambda u^{-\frac{1}{H}}]}\sup_{\tau\in I_k} X(s,\tau)>u}}\\
&\ge
2\mathcal H_{B_H}^{\inf}([0,\lambda a^{\frac{1}{2H}}])
\mathcal H_{B_H}^{\sup}([0,\gamma a^{\frac{1}{2H}}])
\Psi(u) u^{\frac{1}{H}-1}
\frac{\sqrt\pi}{2\gamma\sqrt{b+\varepsilon}}
(1+o(1)),
\end{align*}
which, by the fact that $\varepsilon,\gamma>0$ were arbitrary and $\lim_{S\toi}\frac{1}{S}\mathcal H_{B_H}^{\sup}([0,S])=\mathcal H_{B_H}^{\sup}$,
yield
$$
2\sum_{k=0}^{[\gamma^{-1}u^{\frac{1}{H}}\frac{\log u}{ 2u}]}
\prob{\inf_{s\in[0,\lambda u^{-\frac{1}{H}}]}\sup_{\tau\in I_k} X(s,\tau)>u}
\ge
\mathcal H_{B_H}^{\inf}([0,\lambda a^{\frac{1}{2H}}])
\mathcal H_{B_H}^{\sup}a^{\frac{1}{2H}}
\frac{\sqrt{\pi}}{\sqrt{b}}u^{\frac{1}{H}-1}\Psi(u)
(1+o(1)).
$$
Finally, note that
\begin{align*}
\lefteqn{2\sum_{0\le l<k\le [\gamma^{-1}u^{\frac{1}{H}}\frac{\log u}{ 2u}]}
\mathbb P\Bigg(\inf_{s\in[0,\lambda u^{-\frac{1}{H}}]}\sup_{\tau\in I_k} X(s,\tau)>u,
\inf_{s\in[0,\lambda u^{-\frac{1}{H}}]}\sup_{\tau\in I_l} X(s,\tau)>u\Bigg)}\\
&+
\mathbb P\Bigg(\inf_{s\in[0,\lambda u^{-\frac{1}{H}}]}\sup_{\tau\in I_{-0}} X(s,\tau)>u,
\inf_{s\in[0,\lambda u^{-\frac{1}{H}}]}\sup_{\tau\in I_0} X(s,\tau)>u\Bigg)\\
&\le
2\sum_{0\le l<k\le [\gamma^{-1}u^{\frac{1}{H}}\frac{\log u}{2 u}]}
\mathbb P\Bigg(\sup_{s\in[0,\lambda u^{-\frac{1}{H}}]}\sup_{\tau\in I_k} X(s,\tau)>u,
\sup_{s\in[0,\lambda u^{-\frac{1}{H}}]}\sup_{\tau\in I_l} X(s,\tau)>u\Bigg)\\
&+
\mathbb P\Bigg(\sup_{s\in[0,\lambda u^{-\frac{1}{H}}]}\sup_{\tau\in I_{-0}} X(s,\tau)>u,
\sup_{s\in[0,\lambda u^{-\frac{1}{H}}]}\sup_{\tau\in I_0} X(s,\tau)>u\Bigg).
\end{align*}

It has been shown in \citep[end of the proof of Lemma 3]{Piterbarg01},
that the last expression is of a smaller order than $u^{\frac{1}{H}-1}\Psi(u)$,
which completes the proof of this step.

\subsection{Derivation of the asymptotics}
Recall from \autoref{sec:reduction} that, for any $T>0$,
\[
P(u):=\prob{\inf_{t\in[0,T]}Q_{B_H}(t)>u}=
\prob{\inf_{s\in[0,T A^{\frac{1}{1-H}}(Au^{1-H})^{-\frac{1}{1-H}}]}\sup_{\tau\ge 0}A Z(s,\tau)>Au^{1-H}}.
\]
\autoref{thm:strong} is a simple reformulation of the observations of the
previous steps in terms of the storage process $Q_{B_H}$. We have,
\[
\left[0,T A^{\frac{1}{1-H}}(Au^{1-H})^{-\frac{1}{1-H}}\right]
=
\left[0,\lambda(u)(Au^{1-H})^{-\frac{1}{H}}\right],
\]
where $\lambda(u)= T A^\frac{1}{H} u^{\frac{1-2H}{H}}$. Let $T=T(u)$ be such that $T(u)=o(u^{\frac{2H-1}{H}})$ as $u\toi$. Then, for any $\varepsilon>0$
and all $u$ such that $\lambda(u)\le\varepsilon$,
\[
\left[0,T A^{\frac{1}{1-H}}(Au^{1-H})^{-\frac{1}{1-H}}\right]\subset
\left [0,\varepsilon (Au^{1-H})^{-\frac{1}{H}}\right].
\]
Hence,
\[
P(u)\ge
\prob{\inf_{s\in[0,\varepsilon (Au^{1-H})^{-\frac{1}{H}}]}\sup_{|\tau-\tau_0|\le \log(Au^{1-H})/(Au^{1-H})}A Z(s,\tau)>Au^{1-H}}
\]
and by \eqref{eq.inf1} the last expression is asymptotically bounded below by
\[
\sqrt \pi a^\frac{1}{2H}b^{-\half}\mathcal H_{B_H}^{\sup}\mathcal{ H}^{\inf}_{B_H}\left([0,\varepsilon
a^\frac{1}{2H}]\right)(Au^{1-H})^{\frac{1}{H}-1}\Psi(Au^{1-H}).
\]
Observe that by Fatou's lemma
$\limsup_{\varepsilon\downarrow 0}\mathcal H_{B_H}^{\inf}([0,\varepsilon a^{\frac{1}{2H}}])=1$, which implies the appropriate lower bound for $P(u)$.
Finally, recall from \eqref{eq:asymptotics}, that
\[
\prob{Q_{B_H}(0)>u}
\sim \sqrt\pi a^\frac{1}{2H}b^{-\half}\mathcal H^{\sup}_{B_H}\cdot (Au^{1-H})^{\frac{1-H}{H}}\Psi(Au^{1-H}),\as u,
\]
which is the upper bound for $P(u)$. This completes the proof of \autoref{thm:strong}.

\vspace{1cm}

{\bf Acknowledgement}:  K. D\c{e}bicki was partially supported
by NCN Grant No 2013/09/B/ST1/01778 (2014-2016) and
by the project RARE -318984, a Marie Curie FP7 IRSES Fellowship.

\small\bibliography{biblioteczka}

\end{document}